\documentstyle{amsppt}
\magnification=1200
\UseAMSsymbols
\NoBlackBoxes
\topmatter
\title Total subspaces with long chains of
nowhere norming weak$^*$ sequential closures
\endtitle
\rightheadtext{TOTAL SUBSPACES}
\author M.I.Ostrovskii
\endauthor
\address Mathematical Division, Institute for Low Temperature Physics
and Engineering, 47 Lenin avenue, 310164 Kharkov, UKRAINE \endaddress
\email ostrovskii\%ilt.kharkov.ua\@relay.ussr.eu.net \endemail
\keywords Banach space, Total subspace, Weak$^*$ sequential closure,
Nowhere norming subspace
\endkeywords
\subjclass Primary 46B20\endsubjclass
\abstract If a separable Banach space $X$ is such that for some
nonquasireflexive  Banach  space $Y$  there  exists a surjective
strictly singular operator $T:X\to Y$  then  for  every  countable
ordinal $\alpha $  the dual of $X$  contains  a  subspace  whose weak$^*$
sequential closures of orders less than $\alpha $ are not norming over
any infinite-dimensional subspace of $X$
and whose weak$^*$ sequential closure of order $\alpha +1$ coincides with
$X^*$
\endabstract
%\thanks   \endthanks
\endtopmatter

\document

Let $X$ be a Banach space, $X^*$ be its dual space. The closed unit
ball and the unit sphere of $X$ are denoted  by $B(X)$  and $S(X)$
respectively.  The  term  ``operator''  means  a  bounded   linear
operator.

Let us recall some definitions.  A subspace $M$ of $X^*$ is said to
be {\it total} if for every $0\neq x\in X$ there is an
$f\in M$ such that $f(x)\neq 0$. A
subspace $M$ of $X^*$ is said to be {\it norming over a subspace} $L$ of $X$ if
for some $c>0$ we have
$$
(\forall x\in L)(\sup _{f\in S(M)}|f(x)|\ge c||x||).
$$
\noindent A subspace $M$ of $X^*$ is said to be {\it norming}
if it is  norming  over
$X$. If $M$ is not norming over any infinite dimensional subspace of
$X$ then we shall say that $M$ is {\it nowhere norming}.

The set of all limits  of  weak$^*$  convergent  sequences  in  a
subset $M$ of $X^*$ is called the
{\it weak}$^*$ {\it sequential closure of} $M$ and is
denoted by $M_{(1)}$. If $M$ is a subspace then $M_{(1)}$ is also a subspace.
This subspace need not be closed and all the more need  not  be
weak$^*$ closed \cite{M}. In this connection S.Banach introduced
\cite{B, p.~208,~213}
the weak$^*$ sequential closures (S.Banach used the term
"d\'eriv\'e faible") of other orders, including transfinite ones.
For ordinal $\alpha $ the {\it weak}$^*$
{\it sequential closure of order} $\alpha $ of  a
subset $M$ of $X^*$ is the set
$M_{(\alpha )}=\cup _{\beta <\alpha }(M_{(\beta )})_{(1)}$.

It should be noted that for separable $X$  the  notion  of  the
weak$^*$ sequential closure of order $\alpha $ coincides with the notion of
the derived set of order $\alpha $ considered in \cite{A}, \cite{M1},
\cite{M2}.

For the chain of the weak$^*$ sequential closures we have
$$
M_{(1)}\subset M_{(2)}\subset \ldots
\subset M_{(\alpha )}\subset M_{(\alpha +1)}\subset \ldots.
$$
If we have $M_{(\alpha )}=M_{(\alpha +1)}$ then all subsequent
closures coincide
with $M_{(\alpha )}$. The least ordinal $\alpha $ for which
$M_{(\alpha )}=M_{(\alpha +1)}$ is called
the {\it order} of $M$.

The present paper  deals  with  one  of  the  aspects  of  the
following general problem: how far from the  norming  subspaces
can the total subspaces be and what is the structure  of  Banach
spaces whose duals contain total ``very'' nonnorming subspaces?

There are many works devoted  to  this  problem  (see  \cite{B,  p.~208--215},
papers \cite{A}, \cite{DJ}, \cite{M1}, \cite{M2}, \cite{O1},
\cite{O2} and papers cited  therein).  We
recall only the results which  motivate  us  to  carry  out  the
present research.

1. There is a total subspace $M$ of $l_{\infty }=(l_{1})^*$
such that $M_{(n)}$ is
nowhere norming for all $n\in\Bbb N$ \cite{A}.

2. Let $X$ be a nonquasireflexive separable Banach space. Then,
for every countable ordinal $\alpha $, there  is  a  total  subspace  of
order $\alpha +1$ in $X^*$ \cite{O1}.(V.B.Moscatelli
\cite{M1}, \cite{M2}  obtained  this
result in the case when $\alpha $ is not greater than the first infinite
ordinal. Explicit construction  of \cite{M2} is  useful  in  further
investigation of such subspaces.)

3. Let $X$ be a separable Banach space. Its dual contains a total
nowhere   norming   subspace   if   and   only   if   for   some
nonquasireflexive Banach  space $Y$  there  exists  a  surjective
strictly singular operator $T:X\to Y$ \cite{O2}.

The main result of the present paper is the following.

\proclaim{THEOREM} Let $X$ be a separable Banach space such that for  some
nonquasireflexive Banach  space $Y$  there  exists  a  surjective
strictly singular  operator $T:X\to Y$.  Then  for  every  countable
ordinal $\alpha $ there exists a subspace $M$ of $X^*$  such
that $M_{(\alpha +1)}=X^*$
and for every  ordinal $\beta <\alpha $  the
subspace $M_{(\beta )}\subset X^*$  is  nowhere
norming.
\endproclaim

Let us introduce some notation. For a subset  $A$  of  a  Banach
space $X$, lin$(A)$ and cl$(A)$ are, respectively, the linear span  of
$A$ and the closure of $A$ in the strong topology. By
$w^*-\lim_{m\to \infty }x^*_m$ we
denote the weak$^*$ limit of  the  sequence
$\{x^*_m\}^{\infty }_{m=1}$  in  the  dual
Banach space (if this limit exists). For a subset $A$ of $X^*,\ A^\top$ is
the set $\{x\in X:(\forall x^*\in A)(x^*(x)=0)\}$.

Our sources for Banach space theory are \cite{B}, \cite{S}.

At first we shall prove the theorem in the case when $\alpha $  is  a
nonlimit ordinal. In order to do  this  we  need  the  following
result.
\bigskip
\proclaim{Lemma 1} Let $Y$ be a separable nonquasireflexive Banach  space.
Then for every countable ordinal $\gamma $ there exist a  subspace $N$ of
$Y^*$ and a bounded sequence $\{h_n\}^{\infty }_{n=1}$ in $Y^{**}$ such that:

\noindent$\Bbb A$. If a weak$^*$ convergent sequence
$\{x^*_m\}^{\infty }_{m=1}$ is contained in $N_{(\beta )}$
for some $\beta <\gamma $ and $x^*=w^*-\lim_{m\to \infty }x^*_m$ then
$$
h_n(x^*)=\lim_{m\to \infty }h_n(x^*_m).
\eqno{(1)}$$
$\Bbb B$. There exists a collection
$\{x^*_{n,m}\}^\infty_{n,m=1}$ of vectors in $N_{(\gamma )}$
such that for every $k,n\in\Bbb N$ we have
$$
w^*-\lim_{m\to \infty }x^*_{n,m}=0;
\eqno{(2)}$$
$$
(\forall m\in {\Bbb N})(h_{k}(x^*_{n,m})=\delta _{k,n}).
\eqno{(3)}$$
\endproclaim

At first we shall finish the proof of the theorem in the  case
of nonlimit $\alpha $ with the help of Lemma 1.

We apply Lemma 1 to $\gamma =\alpha -1$ (this ordinal is  correctly  defined
since $\alpha $ is nonlimit). Let $N, \{h_n\}$ and $\{x^*_{n,m}\}$
be as in  Lemma  1.
Let $\{s^*_k\}^{\infty }_{k=1}$ be a normalized sequence in $X^*$ such
that  lin$(\{s^*_k\})$
is a norming subspace of $X^*$. Let $c_1=\sup _n||h_n||$.
Let $\nu _n>0\ (n\in \Bbb N)$  be
such  that $\sum^{\infty }_{n=1}\nu _n\le 1/(2c_1)$.
We  may  assume  without  loss   of
generality that $T$ is a quotient mapping. In this  case $T^*:Y^*\to X^*$
is an isometric embedding. Therefore we may (and shall) identify
$Y^*$ with $T^*(Y^*)\subset X^*$.

Let $R:Y^*\to X^*$ be given by the equality
$$
R(y^*)=y^*+\sum^{\infty }_{n=1}\nu _nh_n(y^*)s^*_n.
\eqno{(4)}$$
It is clear that
$$
(\forall y^*\in Y^*)((1/2)||y^*||\le ||Ry^*||\le (3/2)||y^*||).
\eqno{(5)}$$
Let $M=R(N)$. We prove that for every $\beta \le \gamma $ we have
$$
M_{(\beta )}=R(N_{(\beta )}).
\eqno{(6)}$$
We use the transfinite induction. For $\beta =0$  we have  (6)   by
definition. Let us suppose that (6) is true  for  some
$\beta <\gamma $  and
prove that $M_{(\beta +1)}=R(N_{(\beta +1)})$.
Let $x^*\in M_{(\beta +1)}$, i.e. $x^*=w^*-\lim_{m\to \infty }x^*_m$
for some sequence $\{x^*_m\}^{\infty }_{m=1}$ in $M_{(\beta )}$.
Let $y^*_m\in N_{(\beta )}$  be  such  that
$x^*_m=R(y^*_m)$. Denote $\sup _m||x^*_m||$ by $c_2$. By (5)  we  have
$||y^*_m||\le 2c_2$  for
every $m\in \Bbb N$. Therefore by separability of $X$ we can select a  weak$^*$
convergent    subsequence $\{y^*_{m(i)}\}^{\infty }_{i=1}$
of $\{y^*_m\}^{\infty }_{m=1}$.    Let
$y^*=w^*-\lim_{i\to \infty }y^*_{m(i)}$. It is clear that
$y^*\in N_{(\beta +1)}$. Since $\beta <\gamma $ then by
the assertion $\Bbb A$ of Lemma 1 we have
$$
\lim_{i\to \infty }h_n(y^*_{m(i)})=h_n(y^*)
$$
for every $n\in N$. By the definition of $R$ it follows that
$$
\lim_{i\to \infty }\sum^{\infty }_{n=1}\nu _nh_n(y^*_{m(i)})s^*_n=
\sum^{\infty }_{n=1}\nu _nh_n(y^*)s^*_n,
$$
where the limit is taken in the strong  topology.  Therefore  we
have:
$$
w^*-\lim_{i\to \infty }R(y^*_{m(i)})=R(y^*).
$$
Hence, $x^*=R(y^*)$ and $M_{(\beta +1)}\subset R(N_{(\beta +1)})$.
The inclusion $R(N_{(\beta +1)})\subset M_{(\beta +1)}$
follows immediately  from  (1), (4) and (6).

The case of a limit ordinal $\beta \le \gamma $ is more simple:
$$
M_{(\beta )}=\cup _{\tau<\beta }M_{(\tau)}=
\cup _{\tau<\beta }R(N_{(\tau)})=
R(\cup _{\tau<\beta }N_{(\tau)})=R(N_{(\beta )}).
$$
Therefore  formula  (6)  is  proved.  In  particular,  we   have
$M_{(\gamma )}=R(N_{(\gamma )})$. Let us show that this
equality implies  that $M_{(\gamma )}$ is nowhere norming.

Suppose that it is not the case. Let an  infinite  dimensional
subspace $L$ of $X$ be such that $M_{(\gamma )}$ is norming over $L$.

Recall that if $U$ and $V$ are subspaces of a Banach space $X$  then
the number
$$
\delta (U,V)=\inf \{||u-v||: u\in S(U), v\in V\}
$$
is called the {\it inclination} of $U$ to $V$.

Since $T$ is a strictly singular quotient mapping  then $X$  does
not contain  an  infinite  dimensional  subspace  with  non-zero
inclination to $\ker (T)$. Using the well-known arguments (see  \cite{AGO},
\cite{G} or \cite{R}) we can find a normalized
sequence $\{z_i\}^{\infty }_{i=1}$
in $L$ such
that for some sequence $\{t_i\}^{\infty }_{i=1}$ in
$\ker (T)$ we have $||z_i-t_i||<2^{-i}$ and, furthermore,
$$
(\forall n\in \Bbb N)(\lim_{i\to \infty }s^*_n(t_i)=0).
\eqno{(7)}$$
Let $c>0$ be such that
$$
(\forall x\in L)(\exists f\in S(M_{(\gamma )}))(|f(x)|\ge c||x||).
$$
In particular,
$$
(\forall i\in \Bbb N)(\exists f_i\in S(M_{(\gamma )}))(|f_i(z_i)|\ge c).
$$
Since $M_{(\gamma )}=R(N_{(\gamma )})$  then  we  can
find $y^*_i\in N_{(\gamma )}$  such   that
$f_i=R(y^*_i)$, i.\ e.
$f_i=y^*_i+\sum^{\infty }_{n=1}\nu _nh_n(y^*_i)s^*_n$.
By  (5)  we  have $||y^*_i||\le 2.$
Furthemore, we have
$$
c\le |f_i(z_i)|\le |f_i(z_i-t_i)|+|f_i(t_i)|<
2^{-i}+|\sum^{\infty }_{n=1}\nu _nh_n(y^*_i)s^*_n(t_i)|\le
$$
$$
2^{-i}+c_1||y^*_i||\sum^{\infty }_{n=1}\nu _n|s^*_n(t_i)|.
$$

Using (7) and boundedness of the  sequences $\{s^*_n\}$  and $\{t_i\}$  we
obtain a contradiction. Hence,  the  subspace $M_{(\gamma )}$  is  nowhere
norming. Thus we prove that $M$ satisfies the second assertion  of
the theorem.

It remains to prove that $M_{(\gamma +2)}=X^*$.

Let $\{x^*_{n,m}\}^\infty_{n,m=1}$
be  the  collection  whose  existence   is
asserted in Lemma 1. By  (6)  we  have
$R(x^*_{n,m})\in M_{(\gamma )}$  for  every
$m,n\in \Bbb N$. By (3) we  have $R(x^*_{n,m})=x^*_{n,m}+\nu _ns^*_n$.
Since  the  sequence
$\{x^*_{n,m}\}^{\infty }_{m=1}$ is weak$^*$ null and
$\nu _n\neq 0$  then  we  have $s^*_n\in M_{(\gamma +1)}$  for
every $n\in \Bbb N$,  therefore  lin$(\{s^*_n\})\subset M_{(\gamma +1)}$.
Since  the   subspace
lin$(\{s^*_n\})$ is norming then by \cite{B, p.~213} we have
$M_{(\gamma +2)}=X^*$.  The
theorem is proved in the case when $\alpha $ is nonlimit.

Proof of Lemma 1. By \cite{DJ, Theorem 2} $Y$ contains a bounded  away
from   0   basic   sequence $\{z_n\}^{\infty }_{n=0}$   such
that   the    set
$$\left\{\sum^k_{i=j}z_{i(i+1)/2+j}\right\}^{\infty\ \ \ \infty}_{j=0,k=j}$$
is bounded. We may assume without loss of generality that $||z_i||\le 1$
for every $i\in \Bbb N$. Let $Z=$cl(lin$(\{z_n\}^{\infty }_{n=0}))$.
It is easy  to  see  that
the following claims are true.
\par
1. The space $Z^{**}$ may be identified with the weak$^*$ closure of $Z$
in $Y^{**}$.
\par
2. Every weak$^*$ null sequence in $Z^*$ has a weak$^*$  null  sequence
of extensions to $Y$.
\par
3. If we denote the canonical embedding of $Z$ into $Y$ by $\xi $  then
for every ordinal $\alpha $ and every subspace $N$ of $Z^*$ we shall have
$$
(\xi ^*)^{-1}(N_{(\alpha )})=((\xi ^*)^{-1}N)_{(\alpha )}.
$$
(In this connection see Lemma 1 in \cite{O1}.)
\par
4. If $z^{**}\in Z^{**}$ and $y^*\in Y^*$ then the value $z^{**}(y^*)$ depends
only on
the restriction of $y^*$ to Z.
\par
These claims imply that it is sufficient to prove Lemma 1  with
$Y^*$ and $Y^{**}$ replaced by $Z^*$ and $Z^{**}$ respectively.
\par
Let  us  introduce  some  notation.  We  shall  write $z^{j}_i$  for
$$z_{(j+i-1)(j+i)/2+j}\ (j=0,1,2,\ldots
,\ i\in\Bbb N),$$
biorthogonal  functionals to the system will be
denoted by $\tilde{z}_n(\tilde{z}^{j}_i)$.
By  the  cited  above result of \cite{DJ} we have
$$
\sup _{j,m}||\sum^m_{i=1}z^j_i||=c_1<\infty .
$$
Therefore for every $j=0,1,2,\ldots$ the sequence
$\{\sum^m_{i=1}z^j_i\}^{\infty }_{m=1}$ has  at
least one weak$^*$ limit point in $Z^{**}$. Let us choose one  of  these
limit points and denote it by $f_j$. It is clear that $||f_j||\le c_1$.

We need the following result from \cite{O1}.
\bigskip
\proclaim{Lemma 2} For every vector $g_{0}\in Z^{**}$  of  the  form
$af_j+z^{r}_{s} (a>0,
r\neq j)$, every countable ordinal $\alpha $ and every  infinite  subset
$A\subset \Bbb N$
such that $j,r\not\in A$ there exists a countable subset $\Omega (g_{0},\alpha
,A)$ of $Z^{**}$
such that
\par
1)   For   a    subspace $K(g_{0},\alpha ,A)$    of $Z^*$  defined   by
$K(g_{0},\alpha ,A)=(\Omega (g_{0},\alpha ,A))^\top$ we have $(K(g_{0},\alpha
,A))_{(\alpha )}\subset \ker (g_{0})$.
\par
2) All vectors $h\in \Omega (g_{0},\alpha ,A)$ are of the  form
$h=a(h)f_{j(h)}+z^{r(h)}_{s(h)}$
with $j(h), r(h)\in A\cup \{j,r\}, a(h)>0$  and  for  every $h\neq g_{0}$
from
$\Omega (g_{0},\alpha ,A)$ we have $j(h)\neq r, r(h)\neq $r.
\par
3) If we denote by $Q(b,g_{0},\alpha ,A)$ the intersection of the set
$$
b\tilde{z}^r_s+u,\hbox{ where }u\in \hbox{lin}(\{\tilde{z}^{t}_k\}^{\infty
}_{k=1,t\in A\cup \{j\}})
\eqno{(8)}$$
with $K(g_{0},\alpha ,A)$  then  the  set $(Q(b,g_{0},\alpha ,A))_{(\alpha )}$
contains  all
vectors of the form \rm{(8)} which are in $\ker (g_{0})$.
\endproclaim
\par
Let us introduce the  functionals $h_n\in Z^{**}$  by  the  equalities
$h_n=f_{2n-1}\ (n\in \Bbb N)$ and the vectors $x^*_{n,m}$ by
the equalities $x^*_{n,m}=\tilde{z}^{2n-1}_m$.
It is clear that the vectors $h_n$ and $x^*_{n,m}$ satisfy equalities (2)
and (3).

Let $\{A_n\}^{\infty }_{n=0}$ be a partition of
the set of even natural  numbers
into pairwise disjoint infinite sets. Let $\varepsilon _{n,k}>0\
(n,k\in \Bbb N)$ be such
that $\sum^{\infty }_{n,k=1}\varepsilon _{n,k}<\infty $.
Define  the  family $\{g_{n,k}\}^{\infty }_{n,k=1}$   in   the
following way:
$$
g_{n,k}=z^{2n-1}_k+\varepsilon _{n,k}f_{j(n,k)},
$$
where  the  mapping $j:\Bbb N\times \Bbb N\to \Bbb N$
is   such   that $j(n,k)\in A_n$   and
$j(n,k)\neq j(n,l)$ for $k\neq l$.
\par
Let $\{D_{n,k}\}^{\infty }_{n,k=1}$ be a partition of $A_0$  into  pairwise
disjoint
infinite sets.

The cases of limit and nonlimit $\gamma $ will be treated separately.
\par
Let $\gamma $ be a nonlimit ordinal and let $\Omega (g_{n,k},\gamma
-1,D_{n,k})$  be  the
sets whose existence is asserted in Lemma 2. Let us define $N\subset Z^*$
by $N=(\cup ^{\infty }_{n,k=1}\Omega (g_{n,k},\gamma -1,D_{n,k}))^\top$. Let us
show that $N$ satisfies the
conditions of Lemma 1.
\par
Let $\{x^*_m\}^{\infty }_{m=1}$ be a weak$^*$ convergent sequence in $N_{(\beta
)}(\beta \le \gamma -1)$  and
let $x^*=w^*-\lim_{m\to \infty }x^*_m$. By Lemma 2 we have
$$
(\forall n,k\in \Bbb N)(x^*_m\in \ker (g_{n,k}))
\eqno{(9)}$$
Since the sequence $\{z_n\}^{\infty }_{n=0}$  is  a  basis  of $Z$,  then
their
biorthogonal sequence $\{\tilde{z}_n\}^{\infty }_{n=0}$  is a $w^*$-Schauder
basis  of $Z^*$ \cite{S, p.~155}. (It
means that every vector $z^*\in Z^*$ can be represented
as $z^*=w^*-\lim_{k\to \infty }\sum^k_{n=0}a_n\tilde{z}_n$, where
$a_n=z^*(z_n))$.
\par
Using (9) we can estimate some  of  the  coefficients  of  the
weak$^*$ decompositions of the vectors $x^*_m (m\in \Bbb N)$. Precisely,  if
we
denote $\sup _m||x^*_m||$ by $c_2$, then we obtain
\par
$$
|x^*_m(z^{2n-1}_k)|\le \varepsilon _{n,k}c_1c_2.
$$
Therefore $x^*_m$   can   be   represented   as $u^*_m+v^*_m$,    where
$u^*_m=\sum^{\infty }_{n,k=1}x^*_m(z^{2n-1}_k)\tilde{z}^{2n-1}_k ($we note
that  this  series  converges
unconditionally, therefore we need not  indicate  the  order  of
summability), and
\par
$$
(\forall n,k\in \Bbb N)(v^*_m(z^{2n-1}_k)=0)
\eqno{(10)}$$
Since  the  weak$^*$  convergence  implies   the   coordinatewise
convergence for $w^*$-Schauder bases, then we can represent $x^*$  in
an analogous way, $x^*=u^*+v^*$.
\par
We have
\par
$$
u^*=w^*-\lim_{m\to \infty }u^*_m.
$$
By (10) and by the definition of $h_n$ it follows that
\par
$$
(\forall m,n\in \Bbb N)(h_n(v^*_m)=0);
\eqno{(11)}$$
$$
(\forall n\in \Bbb N)(h_n(v^*)=0).
\eqno{(12)}$$
Since the vectors $\{u^*_m\}^{\infty }_{m=1}$ and $u^*$ are contained in the
strongly
compact set
\par
$$
C=\{\sum^{\infty }_{n,k=1}a_{n,k}\tilde{z}^{2n-1}_k:|a_{n,k}|\le \varepsilon
_{n,k}c_1c_2\},
$$
\noindent then the weak$^*$ convergence  of $\{u^*_m\}$  to $u^*$  implies
the  weak
convergence of $\{u^*_m\}$ to $u^*$.  Therefore  for  every $n\in \Bbb N$  we
have
$\lim_{m\to \infty }h_n(u^*_m)=h_n(u^*)$. From  here  by  (11)  and  (12)  we
obtain
$\lim_{m\to \infty }h_n(x^*_m)=h_n(x^*)$. Thus the assertion $\Bbb A$ of Lemma
1 is proved.
\par
In order to prove the assertion $\Bbb B$ it is  sufficient  to  check
that $\tilde{z}^{2n-1}_k\in N_{(\gamma )}$.  Let
$\tilde{z}(t)=\tilde{z}^{2n-1}_k-\tilde{z}^{j(n,k)}_t/\varepsilon _{n,k} (t\in
\Bbb N)$. It is clear
that $\tilde{z}(t)\in \ker g_{n,k}$ and that $\tilde{z}^{2n-1}_k=w^*-\lim_{t\to
\infty }(\tilde{z}(t))$. Furthermore,
vectors $\tilde{z}(t)$ are of the form  (8)  with $b=1, r=2n-1, s=k$  and
$j=j(n,k)$.    Therefore $\tilde{z}(t)\in (Q(1,g_{n,k},\gamma
-1,D_{n,k}))_{(\gamma -1)}$     and
$\tilde{z}^{2n-1}_k\in (Q(1,g_{n,k},\gamma -1,D_{n,k}))_{(\gamma )}$.
\par
It remains to show that
$$
(\forall r,s\in \Bbb N)(Q(1,g_{n,k},\gamma -1,D_{n,k})\subset (\Omega
(g_{r,s},\gamma -1,D_{r,s}))^\top)
\eqno{(13)}$$
\noindent For $r=n, s=k$ it follows immediately from the  definition  of  $Q$.
Let $(r,s)\neq (n,k)$. Recall that every element  of $\Omega (g_{r,s},\gamma
-1,D_{r,s})$
is a weak$^*$ limit point of linear combinations of $z^{2r-1}_s, z^{j(r,s)}_{t}
(t\in \Bbb N)$ and $z^{p}_{q} (p\in D_{r,s}, q\in \Bbb N)$ and that
$Q(1,g_{n,k},\gamma -1,D_{n,k})$  consists
of  linear   combinations   of $\tilde{z}^{2n-1}_k, \tilde{z}^{j(n,k)}_{t}
(t\in \Bbb N)$  and $\tilde{z}^{p}_{q}
(p\in D_{n,k}, q\in \Bbb N)$.  Our  construction  is  such   that   the   sets
$\{2n-1,j(n,k)\}\cup D_{n,k}$ and $\{2r-1, j(r,s)\}\cup D_{r,s}$ intersect if
and only
if $2n-1=2r-1.$ Since in this case we have $s\neq k$,  then  we  obtain
(13).
\par
Thus we have finished the proof of Lemma 1 in the case when $\gamma $
is a nonlimit ordinal. Let us pass to the case when $\gamma $ is a limit
ordinal. Let $\{\gamma _n\}^{\infty }_{n=1}$ be an increasing sequence of
ordinals,  for
which $\gamma =\lim_{n\to \infty }\gamma_n$  . Let us introduce  the  subspace
$N\subset Z^*$  by  the
equality
$$
N=(\cup ^{\infty }_{n,k=1}\Omega (g_{n,k},\gamma _{n+k},D_{n,k}))^\top.
$$
We shall show that $N$ satisfies all the conditions of Lemma  1.
Let $\{x^*_m\}^{\infty }_{m=1}$ be a weak$^*$ convergent sequence in $N_{(\beta
)}$ with $\beta <\gamma $  and
let $x^*=w^*-\lim_{m\to \infty }x^*_m$. Let $i\in \Bbb N$  be  such  that
$\gamma _{i-1}<\beta \le \gamma _i ($we  let
$\gamma _0=0)$. The definition of the sets $\Omega (g_{n,k},\gamma
_{n+k},D_{n,k})$ implies that
for those pairs $(n,k)$ for which $n+k\ge i$ we have $x^*_m\in \ker (g_{n,k})$,
and,
consequently, we have
$$
|x^*_m(z^{2n-1}_k)|\le \varepsilon _{n,k}c_1c_2.
$$
\noindent In the same time since $||z^{2n-1}_k||\le 1$ then
$$
(\forall n,k\in \Bbb N)(|x^*_m(z^{2n-1}_k)|\le c_2).
$$
\noindent Therefore we may argue in the same way as in the first  part  of
Lemma 1 if we define the set $C$ in the following way:
$$C=\{\sum^{\infty }_{n,k=1}a_{n,k}\tilde{z}^{2n-1}_k:|a_{n,k}|\le \varepsilon
_{n,k}c_1c_2\hbox{ if }n+k\ge i\hbox{ and }$$
$$|a_{n,k}|\le c\hbox{ if } n+k<i\}.$$
The proof of Lemma 1 is complete.
\par
Thus in the case when $\alpha $ is a nonlimit ordinal the proof of the
theorem is finished.
\par
Let us describe the changes which should be made in the  proof
of the theorem in the case when $\alpha $ is a limit ordinal.
\par
Let $\{\alpha _n\}^{\infty }_{n=1}$ be an increasing  sequence  of  nonlimit
ordinals
such that $\alpha =\lim_{n\to \infty }\alpha _n$. Instead  of  Lemma  1  we
shall  use  the
following result.
\proclaim{Lemma 3} For every separable nonquasireflexive Banach space $Y$
there exists a subspace $N$ of $Y^*$ and a bounded sequence $\{h_n\}^{\infty
}_{n=1}$
in $Y^{**}$ such that:

\noindent$\Bbb A^{\hbox{new}}$. If a weak$^*$ convergent sequence
$\{x^*_m\}^{\infty }_{m=1}$ is  contained  in
$N_{(\beta )}$ for some $\beta <\alpha $ and $x^*=w^*-\lim_{m\to \infty }x^*_m$
then we have
$$
h_n(x^*)=\lim_{m\to \infty }h_n(x^*_m)
$$
for those $n$ for which $\beta <\alpha _n$.

\noindent$\Bbb B^{\hbox{new}}$. For every $n\in \Bbb N$
there exists a sequence $\{x^*_{n,m}\}^{\infty }_{m=1}$ in
$N_{(\alpha_n)}$ such that the conditions \rm{(2)} and \rm{(3)} are satisfied.
\endproclaim
\par
At first we finish the proof of the theorem with the  help  of
Lemma 3. Let $\{s^*_k\}^{\infty }_{k=1}$ and $\{\nu _n\}^{\infty }_{n=1}$
be the same  as  in  the  first part of the theorem.
Let $R:Y^*\to X^*$ be defined by equality (4). For
every ordinal $\beta <\alpha $ we denote by $W(\beta )$ the set of natural
numbers $n$
for which $\alpha _n<\beta $. It is clear that $W(\beta )$ is a finite set.
\par
Let us show that for every ordinal $\beta <\alpha $ we have
$$
M_{(\beta )}=\hbox{lin}(R(N_{(\beta )})\cup \{s^*_k\}_{k\in W(\beta )}).
\eqno{(14)}$$
We shall prove this with the aid of transfinite induction. For
$\beta =0$ (14) follows immediately from the definition. Let us suppose
that (14)  is  valid  for  some $\beta <\alpha $  and  prove  the  analogous
equality for $\beta +1.$ Let $x^*\in M_{(\beta +1)}$,  i.e. $x^*=w^*-\lim_{m\to
\infty }x^*_m$,  where
$x^*_m\in M_{(\beta )}$=lin$(R(N_{(\beta )})\cup \{s^*_k\}_{k\in W(\beta )})$.
It is clear that the last space
is a subspace of cl$(R(N_{(\beta )}))\oplus F$, where $F$  is  some  subspace
of
lin$(\{s^*_k\}_{k\in W(\beta )})$.  Therefore,  we  may  write
$x^*_m=u^*_m+v^*_m$,   where
$u^*_m\in$cl$(R(N_{(\beta )}))$ and $v^*_m\in $F. It is clear that the
sequences $\{u^*_m\}^{\infty }_{m=1}$
and $\{v^*_m\}^{\infty }_{m=1}$ are bounded and that we may suppose without
loss  of
generality that $u^*_m\in R(N_{(\beta )})$. Let $u^*_m=Ry^*_m$. By (5) the
sequence $\{y^*_m\}$
is also bounded. So we can find a sequence $\{m(i)\}^{\infty }_{i=1}$ of
natural
numbers such that the sequences $\{v^*_{m(i)}\}^{\infty }_{i=1}$  and
$\{y^*_{m(i)}\}^{\infty }_{i=1}$  are
weak$^*$ convergent. Let $v^*$ and $y^*$ be corresponding  weak$^*$  limits.
By Lemma 3 we have $h_n(y^*)=\lim_{i\to \infty }h_n(y^*_{m(i)})$ for all $n$
for  which
$\alpha _n>\beta $. Therefore,
\par
$$
\lim_{i\to \infty }\sum^{}_{n\not\in W(\beta +1)}\nu
_nh_n(y^*_{m(i)})s^*_n=\sum^{}_{n\not\in W(\beta +1)}\nu _nh_n(y^*)s^*_n,
$$
\noindent where the limit is taken in the strong topology.
\par
Since the set $W(\beta +1)$ is finite then we  may  without  loss  of
generality assume that the sequence $\{\sum^{}_{n\in W(\beta +1)}\nu
_nh_n(y^*_{m(i)})s^*_n\}^{\infty }_{i=1}$
is strongly convergent.  Let $\sum^{}_{n\in W(\beta +1)}a_ns^*_n$  be  its
limit.  It
follows that
\par
$$
x^*=R(y^*)-\sum^{}_{n\in W(\beta +1)}\nu _nh_n(y^*)s^*_n+\sum^{}_{n\in W(\beta
+1)}a_ns^*_n+v^*\in
$$
$$
\hbox{lin}(R(N_{(\beta +1)})\cup \{s^*_n\}_{n\in W(\beta +1)}).
$$
Therefore $M_{(\beta +1)}\subset $lin$(R(N_{(\beta +1)})\cup \{s^*_n\}_{n\in
W(\beta +1)})$.  The   converse
inclusion  follows  immediately  from $\Bbb B^{\hbox{new}}$  and  the
induction
hypothesis.
\par
If $\beta $ is a limit ordinal and (14) is proved  for  all  ordinals
less than $\beta $ then (14) follows  immediately  from  the  following
fact: $W(\beta )=W(\tau)$ for some $\tau<\beta $. So (14) is proved for all
ordinals
which are less than $\alpha $.
\par
It is not hard to check that the linear span of the union of a
finite-dimensional and a nowhere norming  subspaces  is  nowhere
norming. Therefore by (14) and the arguments of the  first  part
of the theorem $M_{(\beta )}$ is a nowhere norming subspace for every $\beta
<\alpha $.
\par
The proof of the equality $M_{(\alpha +1)}=X^*$ is  the  same  as  in  the
first part of the theorem.
\par
Proof of Lemma 3. We repeat the  arguments  of  the  proof  of
Lemma 1 up to the passage where $A_{0}$ was  presented  in  the  form
$A_{0}=\cup ^{\infty }_{n,k=1}D_{n,k}$. Now we continue in the following way.
Let $\{\alpha _n\}^{\infty }_{n=1}$
be the  increasing  sequence  of  nonlimit  ordinals  introduced
above. Let $\Omega (g_{n,k},\alpha _n-1,D_{n,k})$ be  the  sets  whose
existence  is
guaranteed by Lemma 2. Let us introduce the subspace $N\subset Z^*$ by the
equality
\par
$$
N=(\cup ^{\infty }_{n,k=1}\Omega (g_{n,k},\alpha _n-1,D_{n,k}))^\top.
$$
The assertion $\Bbb A^{\hbox{new}}$ of Lemma 3 is proved in the same way as the
assertion $\Bbb A$ of Lemma 1. The only distinction  is  that  we  have
relation (9) not for all natural $n$ but only for $n\in \Bbb N\backslash
W(\beta +1)$. This
do not prevent to finish the  proof  since  the  assertion $\Bbb
A^{\hbox{new}}$
concerns only those $h_n$ for which $n\in \Bbb N\backslash W(\beta +1)$.
\par
The proof of the assertion of $\Bbb B^{\hbox{new}}$ is the same as the proof of
the assertion $\Bbb B$ of Lemma 1.

\enddocument